\documentclass[11pt]{article}

\usepackage{graphicx}
\usepackage{amsmath}
\usepackage{amsfonts}
\usepackage{fancyhdr} 
\usepackage{mathrsfs}
\usepackage{amssymb}
\usepackage{wrapfig}
\usepackage{bbold}
\usepackage{amsthm}

\newtheorem{thm}{Theorem}

\theoremstyle{definition}
\newtheorem{alg}{Algorithm}

\numberwithin{equation}{section}
\numberwithin{defn}{section}

\newcommand{\real}{\mathbb{R}}

\newcommand{\integer}{\mathbb{Z}}

\DeclareMathOperator{\diag}{diag}

\providecommand{\norm}[1]{\lVert#1\rVert}

\newenvironment{pf}{\noindent {\em Proof}.\ \ }{\hspace*{\fill}\rule{.5ex}{1.4ex}\,}

\begin{document}

\title{A scattering-based algorithm for wave propagation in one dimension}
\author{Peter C.~Gibson\footnote{Dept.~of Mathematics \& Statistics, York University, 4700 Keele St., Toronto, Ontario, Canada, M3J~1P3, $\mathtt{pcgibson@yorku.ca}$}}
\date{August 25, 2017}
\makeatletter
\let\newtitle\@title
\let\newauthor\@author
\let\newdate\@date
\makeatother

\maketitle

\begin{center}MSC 35L05, 65M; Keywords: one dimensional wave equation, numerical methods, scattering\end{center}

\begin{abstract}
We present an explicit numerical scheme to solve the variable coefficient wave equation in one space dimension with minimal restrictions on the coefficient and initial data.  
\end{abstract}


\section{The algorithm\label{sec-scattering-algorithm}}

A range of physical and biological applications involve the one dimensional wave equation\footnote{Note that the equation $u_{tt}-c^2u_{zz}=0$ (with $\zeta=1/c$) and the coupled system $\rho u_t+p_z=0$ and $\frac{1}{K}p_t+u_z=0$  (with $c=\sqrt{K/\rho}$, $\zeta=\sqrt{K\rho}$) transform to (\ref{wave-eqn}) under the change of variables $x=\int_0^z\frac{1}{c(s)}\,ds$.}
\begin{subequations}\label{wave}
\begin{gather}
u_{tt}-\textstyle\frac{1}{\zeta}(\zeta u_x)_{x}=0\label{wave-eqn}\\
u(x,0)=f(x),\qquad u_t(x,0)=g(x)\label{ic}
\end{gather}
\end{subequations}
where: there exist $x_-<x_+$ such that $\zeta$ takes constant values $\zeta_-$ and $\zeta_+$ on the respective intervals $(-\infty,x_-)$ and $(x_+,\infty)$; and $0<c<\zeta<C$ for some $c<C\in\real$.  Applications include imaging of layered media such as seismic imaging and the acoustic imaging of laminated structures, microwave imaging of skin tissue, and modelling the human vocal tract and cochlea; see \cite{FoGaPaSo:2007},\cite{BlCoSt:2001}, \cite{WiFeWe:2006}, \cite{vdDAs:2008}, \cite{ZwSh:1995}.  Certain of these applications have a long history, yet there seems to be no clear consensus on how best to compute solutions to (\ref{wave}) in the general setting where $\zeta$ has complicated structure on $[x_-,x_+]$ such as rapid oscillations or discontinuities, or where the initial data is singular (e.g., a Dirac comb), or both (for background see \cite[Ch.~10]{Le:2007}, \cite[Ch.~5]{As:2008}, \cite[Ch.~9]{Gu:2008}).  Recent theoretical insights into the case where $\zeta$ is piecewise constant offer a fresh perspective from which to consider the numerical analysis of  (\ref{wave}) for general $\zeta$ (\cite{Gi:SIAP2014}, \cite{Gi:JFAA2016}). 

The purpose of this note is to present a simple, explicit numerical scheme to solve (\ref{wave}) on a region $(x,t)\in[a,b]\times[0,T]\subset\real^2$ where $[x_-,x_+]\subset[a,b]$.  The scheme is based on a natural idea from scattering theory, but with crucial differences from established methods in both the details of the implementation and the interpretation.  These differences yield an order of magnitude reduction in computational expense coupled with remarkable versatility.  

\pagestyle{fancyplain}
\begin{alg}\label{alg-main}\ \\ 
\emph{Preliminary step.}  Convert the initial data $f(x)$ and $g(x)$ from (\ref{ic}) to a pair of functions
\begin{equation}\label{initial}
\alpha(x)=\frac{1}{2}\left(f(x)-\frac{1}{\zeta(x)}\int_{-\infty}^x\zeta(s)g(s)\,ds\right),\;\beta(x)=\frac{1}{2}\left(f(x)+\frac{1}{\zeta(x)}\int_{-\infty}^x\zeta(s)g(s)\,ds\right).
\end{equation}
\hspace*{\fill}\begin{tabular}{|l|l|}
\hline
\multicolumn{2}{|c|}{\textbf{The scattering algorithm}}\\
\hline
\emph{Input}&\rule{0pt}{12pt}$\zeta,\alpha,\beta,a,b,n,T$\\[5pt]
\emph{Spatial grid}&$\Delta=(b-a)/n,\quad x_j=a+(j-1)\Delta\qquad(1\leq j\leq n+1)$\\[5pt]
\emph{Temporal grid}&$m=\lfloor T/\Delta\rfloor,\quad t_k=k\Delta\qquad(0\leq k\leq m)$\\[5pt]
\emph{Weights}& $r_j=\bigl( \zeta(x_j)-\zeta(x_{j+1})\bigr)/\bigl( \zeta(x_j)+\zeta(x_{j+1})\bigr)\qquad(1\leq j\leq n)$\\[5pt]
\emph{Variables}&$v=\bigl(v^k_i\bigr)\in\real^{(m+1)\times (2n+2)},\quad \tilde u=\bigl(\tilde{u}^k_j\bigr)\in\real^{(m+1)\times (n+1)}$\\[5pt]
\emph{Initialization}&$v^0_{2j-1}=\alpha(x_j),\quad v^0_{2j}=\beta(x_j),\quad \tilde{u}^0_j=v^0_{2j-1}+v^0_{2j}\qquad(1\leq j\leq n+1)$\\[5pt]
&$v^k_1=\alpha(x_1-t_k),\quad v^k_{2n+2}=\beta(x_{n+1}+t_k)\qquad(1\leq k\leq m)$\\[5pt]
\emph{Iterative step}& Given $v^k$ $(0\leq k\leq m-1)$, set\\[5pt]
&$v^{k+1}_{2j-1}=(1+r_{j-1})v^k_{2j-3}-r_{j-1}v^k_{2j}\qquad(2\leq j\leq n+1)$\\[5pt]
&$v^{k+1}_{2j}=r_jv^k_{2j-1}+(1-r_j)v^k_{2j+2}\qquad(1\leq j\leq n)$\\[5pt]
&$\tilde{u}^{k+1}_j=v^{k+1}_{2j-1}+v^{k+1}_{2j}\qquad(1\leq j\leq n+1)$\\[5pt]
\emph{Output}&$\tilde{u}=\bigl(\tilde{u}^k_j\bigr)\in\real^{(m+1)\times (n+1)}$\\[5pt]
\hline
\end{tabular}\hspace*{\fill}
\end{alg}

\section{Discussion\label{sec-remarks}}

\begin{figure}[p]
\hspace*{\fill}
\parbox{6in}{
\includegraphics[width=1.5in]{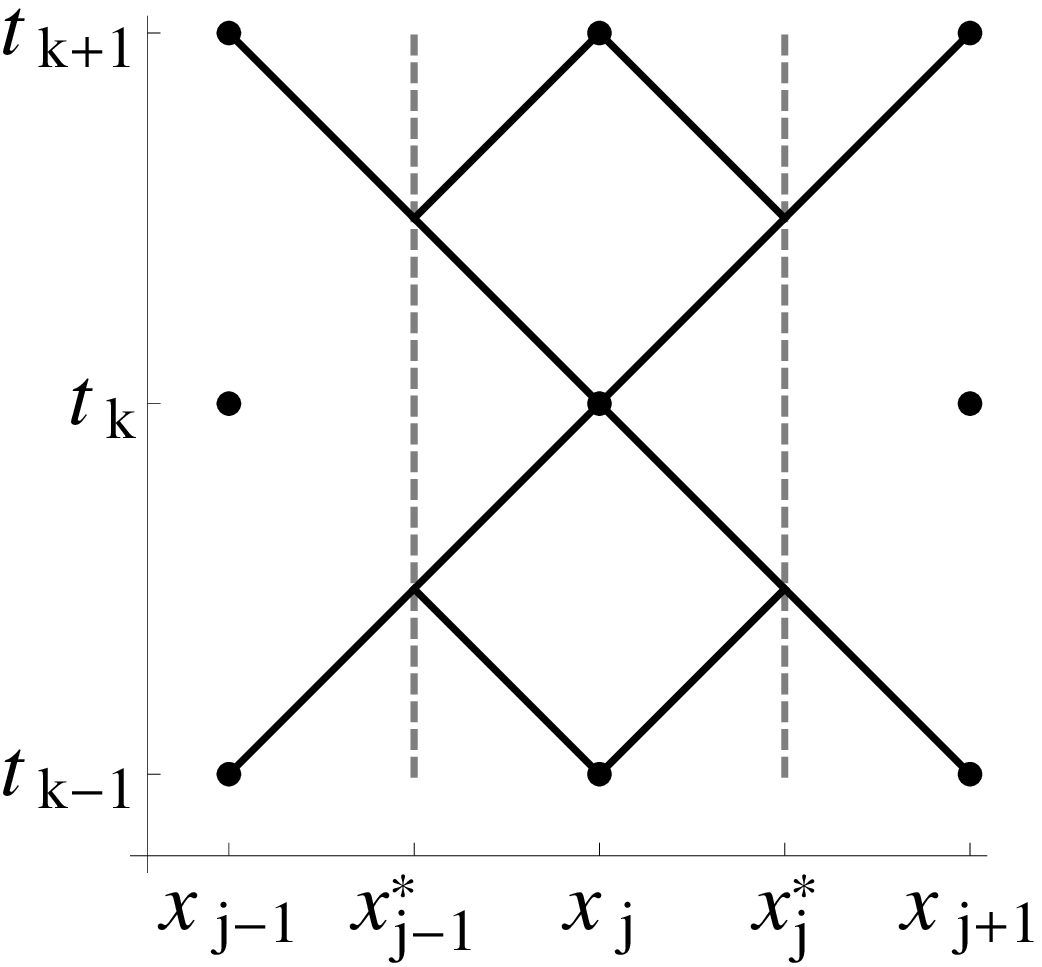}\qquad\includegraphics[width=1.375in,clip=true,trim=0pt -15pt 0pt 0pt]{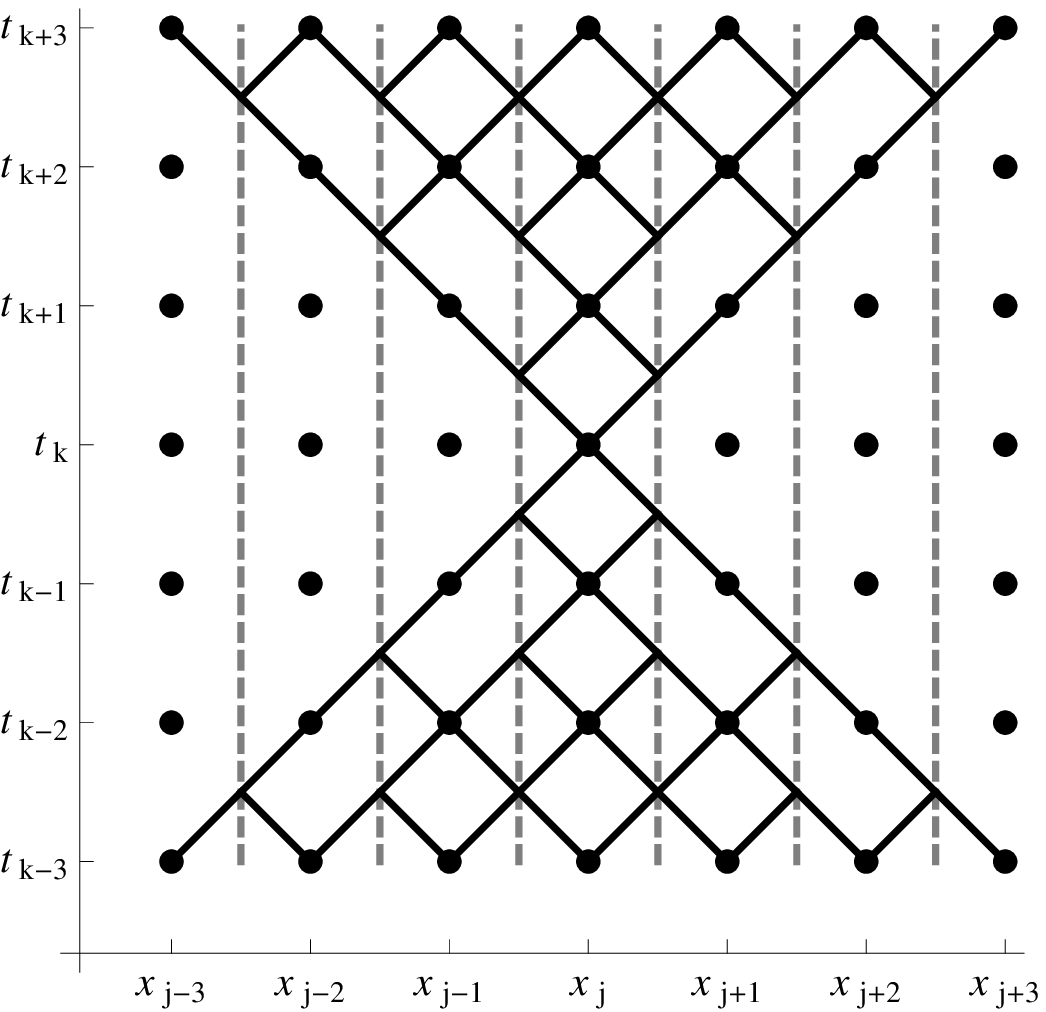}\qquad\includegraphics[width=2in,clip=true,trim=0pt -20pt 0pt 0pt]{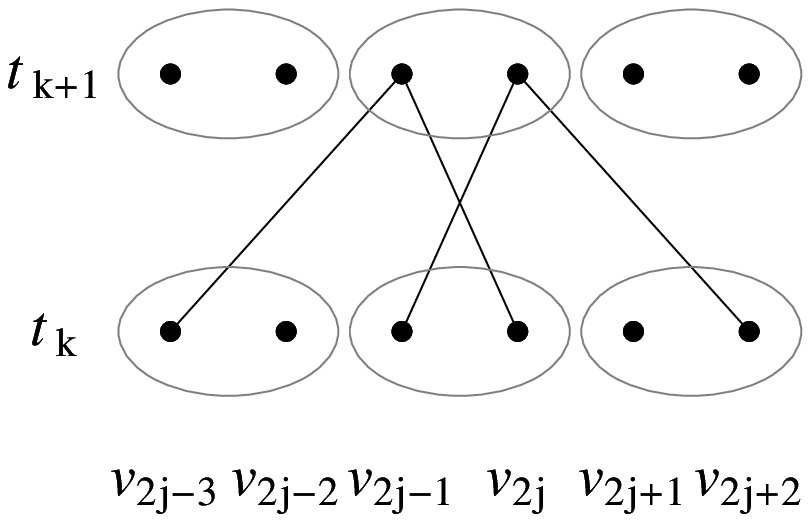}
\caption{Forward and backward characteristics through a grid point $(x_j,t_k)$ in the presence of discontinuities. Left: one time step, center: several time steps. Right: the stencil for Algorithm~\ref{alg-main} with respect to the values $v^k_i$, pairs of which correspond to spatial grid points.}\label{fig-cones}
}
\hspace*{\fill}
\end{figure}

Rather than discretize derivatives, the idea behind Algorithm~\ref{alg-main} is to replace $\zeta$ with a step function and solve the resulting equation exactly.  This echoes the piecewise constant approximations of Goupillaud and Godunov and their many subsequent manifestations (cf.~\cite{BeGoWa:1958}, \cite{Go:1960}, \cite{BuBu:1983}, \cite{FoLe:1999}, \cite[Ch.~3]{FoGaPaSo:2007});
however, Algorithm~\ref{alg-main} differs markedly from previously established methods.  In detail, let $x_j^\ast=(x_j+x_{j+1})/2$ $(1\leq j\leq n)$ denote the midpoints of consecutive spatial grid points, write 
$
P=\{x_1,x_1^\ast,x_2,x_2^\ast,\ldots,x_n,x_n^\ast,x_{n+1}\},
$
and set 
\begin{equation}\label{zeta-P}
\zeta^P(x)=\left\{\begin{array}{ll}
\zeta_-&\mbox{ if }x<x_1^\ast\\
\zeta(x_j)&\mbox{ if } x_{j-1}^\ast\leq x<x_j^\ast\quad(2\leq j\leq n)\\
\zeta_+&\mbox{ if }x_n^\ast\leq x
\end{array}\right..
\end{equation}
Then $\zeta^P$ is a step function having equally spaced jump points $x_1^\ast,\ldots,x_n^\ast$, and it coincides with $\zeta$ on the spatial grid.  
Replacing $\zeta$ by $\zeta^P$ in (\ref{wave}) yields 
\begin{subequations}\label{approximate-wave}
\begin{gather}
u_{tt}-\textstyle\frac{1}{\zeta^P}(\zeta^P u_x)_{x}=0\label{approx-wave}\\
u(x,0)=f(x),\qquad u_t(x,0)=g(x)\label{approx-ic},
\end{gather}
\end{subequations}
the first equation of which reduces to $u_{tt}-u_{xx}=0$ on intervals of the form $(-\infty,x_1^\ast)$, $(x_{j-1}^\ast,x_j^\ast)$ or $(x_{n}^\ast,\infty)$, where $\zeta^P$ is constant. Its solution within each such interval is therefore a pair of left and right travelling waves.  One-sided limits of $u$ and $\zeta^Pu_x$ are required to agree at each jump point $x_j^\ast$
\[
u(x_j^\ast-,t)=u(x_j^\ast+,t)\quad\mbox{ and }\quad\zeta^P(x_j)u_x(x_j^\ast-,t)=\zeta^P(x_{j+1})u_x(x_j^\ast+,t).
\]
By a standard computation (see \cite[Ch.~3]{FoGaPaSo:2007}), this yields reflection and transmission factors: a right-travelling wave in $(x_{j-1}^\ast,x_{j}^\ast)$ gets partially reflected back off $x_j^\ast$, with reflection factor $r_j$ as defined in Algorithm~\ref{alg-main}, and partially transmitted into $(x_j^\ast,x_{j+1}^\ast)$, with transmission factor $1+r_j$.  In particular, characteristic lines for the equation (\ref{approx-wave}) bifurcate at jump points of $\zeta^P$. The corresponding domains of dependence and influence for a grid point $(x_j,t_k)$ are depicted in Figure~\ref{fig-cones}.   In Algorithm~\ref{alg-main} the pair $v^k_{2j-1}, v^k_{2j}$ encodes the respective amplitudes of the right and left moving parts of the solution to (\ref{approximate-wave}) at the grid point $(x_j,t_k)$, and the iterative step accounts for the reweighted reflected and transmitted waves comprising the domain of dependence of $(x_j,t_{k+1})$ at time $t_k$.   In this way the algorithm computes the \emph{exact} solution to (\ref{approximate-wave}) at the grid points.   

The matrix formulation of the dependence of (row vector) $v^{k+1}$ on $v^k$ as expressed in the scattering algorithm is
\begin{equation}\label{matrix-formulation}
v^{k+1}=v^kA_P +(v^{k+1}_1,0,\ldots,0,v^{k+1}_{2n+2}),
\end{equation}
where $A_P$ is a $(2n+2)\times(2n+2)$ pentadiagonal matrix having $4n$ non-zero entries.  For example, if $n=4$,
\begin{equation}\label{AP}
A_P=
\left(
\begin{array}{cccccccccc}
 0 & r_1 & r_1+1 & 0 & 0 & 0 & 0 & 0 & 0 & 0 \\
 0 & 0 & 0 & 0 & 0 & 0 & 0 & 0 & 0 & 0 \\
 0 & 0 & 0 & r_2 & r_2+1 & 0 & 0 & 0 & 0 & 0 \\
 0 & 1-r_1 & -r_1 & 0 & 0 & 0 & 0 & 0 & 0 & 0 \\
 0 & 0 & 0 & 0 & 0 & r_3 & r_3+1 & 0 & 0 & 0 \\
 0 & 0 & 0 & 1-r_2 & -r_2 & 0 & 0 & 0 & 0 & 0 \\
 0 & 0 & 0 & 0 & 0 & 0 & 0 & r_4 & r_4+1 & 0 \\
 0 & 0 & 0 & 0 & 0 & 1-r_3 & -r_3 & 0 & 0 & 0 \\
 0 & 0 & 0 & 0 & 0 & 0 & 0 & 0 & 0 & 0 \\
 0 & 0 & 0 & 0 & 0 & 0 & 0 & 1-r_4 & -r_4 & 0
\end{array}
\right).
\end{equation}
Solutions to (\ref{approximate-wave}) outlined in, e.g., \cite{BuBu:1983} and \cite{FoGaPaSo:2007} involve upper or lower triangular systems requiring approximately $Tn^3/(b-a)$ flops, while the sparse structure of $A_P$ reduces this to $7Tn^2/(b-a)$ flops.  Thus Algorithm~\ref{alg-main} is qualitatively faster than previously established methods, allowing computation on a vastly finer spatial grid.  

The following facts justify Algorithm~\ref{alg-main} theoretically.
\begin{thm}\label{thm-one}
The eigenvalues of $A_P$ lie properly inside the unit circle, implying stability of the scattering algorithm.  \end{thm}
\begin{pf}
Let $U_P$ denote the $(2n+2)\times(2n+2)$ pentadiagonal matrix
\begin{equation}\label{UP}
\begin{split}
&U_P=\\
&\left(
\begin{array}{cccccccccc}
 0 & r_1 & \sqrt{1-r_1^2} & & & & & & & \\
 1 & 0 & 0 & 0 & & & & &\mbox{\Large $0$} & \\
 0 & 0 & 0 & r_2 & \sqrt{1-r_2^2} & & & & & \\
 & \sqrt{1-r_1^2} & -r_1 & 0 & 0 & 0 & & & & \\
 & & 0 & 0 & 0 & \ddots & \ddots & & & \\
 & & &\sqrt{1-r_2^2} & -r_2 & \ddots & 0 & 0 & & \\
 & & & & \ddots & \ddots & 0 & r_n & \sqrt{1-r_n^2} & \\
 & &  & & & \sqrt{1-r_{n-1}^2} & -r_{n-1} & 0 & 0 & 0 \\
 &\mbox{\Large $0$}& & & & & 0 & 0 & 0 & 1 \\
 & & & & & & & \sqrt{1-r_n^2} & -r_n & 0
\end{array}
\right).
\end{split}
\end{equation}
$U_P$ is evidently unitary since its columns (and rows) are pairwise orthogonal.  Set 
\[
\sigma_j=\prod_{k=j}^n\sqrt{\frac{1+r_k}{1-r_k}}=\sqrt{\frac{\zeta(x_j)}{\zeta(x_{n+1})}}\qquad(1\leq j\leq n),
\]
and let $D_P$ denote the $(2n+2)\times(2n+2)$ diagonal matrix 
\begin{equation}\label{DP}
D_P=\diag(\sigma_1,\sigma_1,\sigma_2,\sigma_2,\ldots,\sigma_n,\sigma_n,1,1).
\end{equation}
Let $J$ denote the $(2n+2)\times(2n+2)$ diagonal projection of rank $2n$ obtained by zeroing out the first and last diagonal entries of the identity matrix,
\[
J=\diag(0,1,1,\ldots,1,0).
\]
Then by direct computation,
\[
A_P=D_P(U_PJ)D_P^{-1}.
\]
Thus the spectrum of $A_P$ is the same as that of $U_PJ$.  Note that for any column vector $v$, $\norm{Jv}\leq\norm{v}$, with equality only if $Jv=v$.  If $U_PJv=\lambda v$ for some scalar $\lambda$, then, since $U_P$ is unitary,  
\[
|\lambda|\norm{v}=\norm{U_PJv}=\norm{Jv}\leq\norm{v}.
\]
Hence the spectrum of $U_PJ$ belongs to the closed unit disk.  If $|\lambda|=1$, then $Jv=v$, as noted above, and hence $U_Pv=\lambda v$.  Given that $|\lambda|=1$, this forces $v=0$, as follows. Note that $v_1=v_{2n+2}=0$ since $Jv=v$.  Since $U_Pv=\lambda v$, the structure of the second row of $U_P$ in (\ref{UP}) implies that $v_1=\lambda v_2$, so $v_2=0$.  Looking at the first row of (\ref{UP}), $r_1v_2+\sqrt{1-r_1^2}v_3=\lambda v_1$, forcing $v_3=0$.  In general, if $v_{2k}=v_{2k+1}=0$ for some $1\leq k\leq n-1$, then inspection of row $2k+2$ of $U_P$, followed by $2k+1$ yields that $v_{2k+2}=v_{2k+3}=0$.  Since $v_{2n+2}=0$ this proves $v=0$ as claimed.  Thus if $|\lambda|=1$ then $\lambda$ is not an eigenvalue, so the spectrum of $U_PJ$ and hence that of $A_P$, is contained in the interior of the unit disk. 
\end{pf}
\begin{thm}\label{thm-two}
Let $u^P$ denote the solution to (\ref{approximate-wave}), $(x_j,t_k)$ an arbitrary grid point as defined in Algorithm~\ref{alg-main}, and $\tilde u$ the output of the algorithm. 
\begin{enumerate}
\item If $\alpha$ and $\beta$ are regular functions (as opposed to distributions), then $u^P(x_j,t_k)=\tilde{u}^k_j$.  
\item If $\alpha(x)=\sum_{j\in\integer} c_j\delta(x-a-j\Delta)$ and $\beta(x)=\sum_{j\in\integer} d_j\delta(x-a-j\Delta)$ are Dirac combs supported on an extension of the spatial grid, and values of $\alpha$ and $\beta$ are replaced in Algorithm~\ref{alg-main} by the corresponding coefficients $c_j$ and $d_j$, then $u^P(x_j,t_k)=\tilde{u}^k_j\;\delta(x-a-j\Delta)$.  
\end{enumerate}
\end{thm}
\begin{pf}
Let $\alpha_k$ and $\beta_k$ denote the respective right and left-moving components of $u^P$ at time $t_k$ analogous to (\ref{initial}),
\[
\begin{split}
\alpha_k(x)&=\frac{1}{2}\left(u^P(x,t_k)-\frac{1}{\zeta(x)}\int_{-\infty}^x\zeta(s)u^P_t(s,t_k)\,ds\right)\mbox{ and }\\\beta_k(x)&=\frac{1}{2}\left(u^P(x,t_k)+\frac{1}{\zeta(x)}\int_{-\infty}^x\zeta(s)u^P_t(s,t_k)\,ds\right).
\end{split}
\]
The standard laws of reflection and transmission described above imply that for $x\in(x_{j-1}^\ast,x_j^\ast)$, 
\[
\begin{split}
\alpha_{k+1}(x)&=(1+r_{j-1})\alpha_k(x-\Delta)-r_{j-1}\beta_k(x_j^\ast-(x-x_{j-1}^\ast))\quad\mbox{ and }\\
\beta_{k+1}(x)&=r_j\alpha_k(x_j^\ast-(x-x_{j-1}^\ast))+(1-r_j)\beta_k(x+\Delta).
\end{split}
\]
See Figure~\ref{fig-cones}.  These formulas can be given a valid interpretation for distributions (evaluated on appropriately chosen test functions) as well as ordinary functions.  
Specializing to the point $x=x_j\in(x_{j-1}^\ast,x_j^\ast)$ yields the formulas
\begin{subequations}\label{point}
\begin{align}
\alpha_{k+1}(x_j)&=(1+r_{j-1})\alpha_k(x_{j-1})-r_{j-1}\beta_k(x_j)\quad\mbox{ and }\label{alpha-point}\\
\beta_{k+1}(x_j)&=r_j\alpha_k(x_j)+(1-r_j)\beta_k(x+\Delta).\label{beta-point}
\end{align}
\end{subequations}
The latter formulas make sense for ordinary functions $\alpha_k,\beta_k$. But distributions do not in general have meaningful restrictions to a single point---unless they are supported at a single point, such as is the case for Dirac functions. Thus formulas (\ref{point}), which are precisely the iterative step of the scattering algorithm, apply both to ordinary functions and Dirac distributions, proving that the algorithm coincides with the exact solution to (\ref{wave}) in these cases. 
\end{pf}
\begin{thm}\label{thm-three}
If there is a sequence of partitions $\{P\}$ such that $\zeta^P\rightarrow\zeta$ uniformly, then $u^P\rightarrow u$ on $[a,b]\times[0,T]$ in an appropriate sense depending on the smoothness of the initial data.  
\end{thm}
The convergence $u^P\rightarrow u$ is uniform if initial data are smooth; for certain distributional initial data the convergence may be in the Sobolev space $H^{-1}$.  Proofs and full technical details will be presented elsewhere.  In essence, Theorem~\ref{thm-three} follows from abstract functional analysis such as \cite[Ch.~3.8]{LiMa:1972} or \cite[Ch.~10]{Br:2011}.  See \cite[Theorem~3.5]{KiRi:2016} for a particular version, proved by way of the Hille-Yosida theorem. The heuristic upshot of Theorem~\ref{thm-three} is that uniform convergence of the coefficient $u^P$ drives convergence of the solution.  

\begin{figure}[p]
\hspace*{\fill}
\parbox{6in}{
\includegraphics[width=2in]{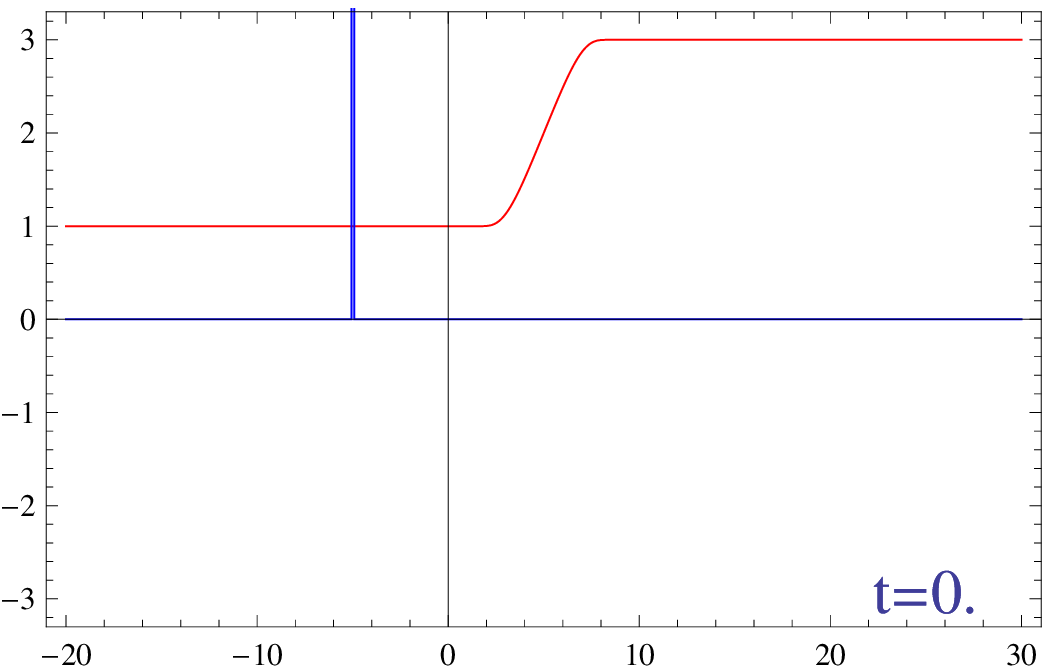}\includegraphics[width=2in]{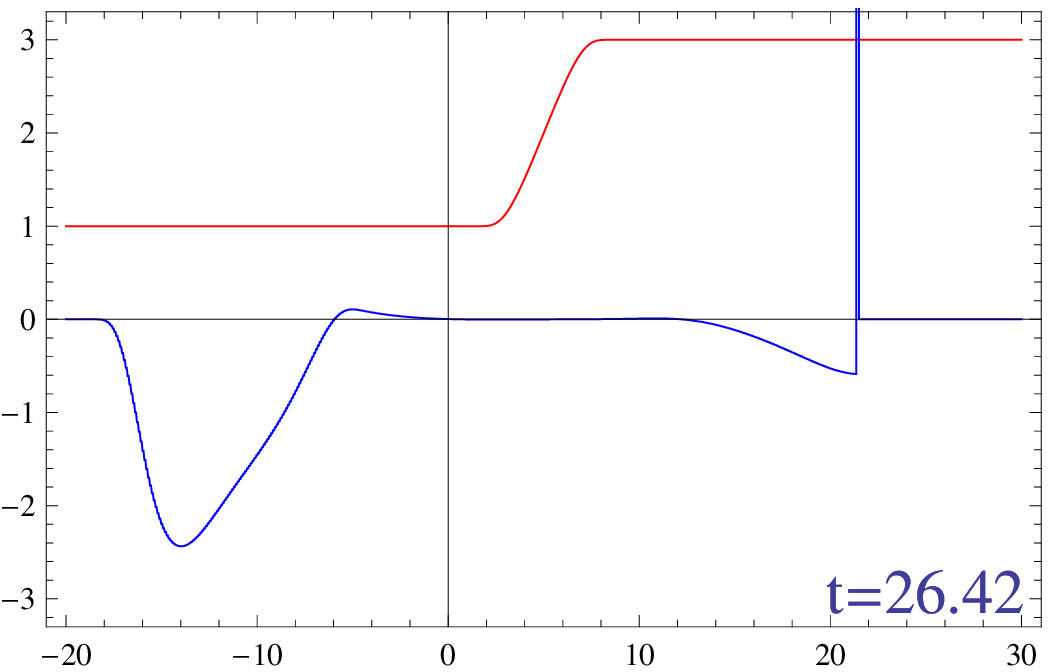}\quad\includegraphics[width=1.5in, clip=true,trim=0pt -25pt 0pt .0pt]{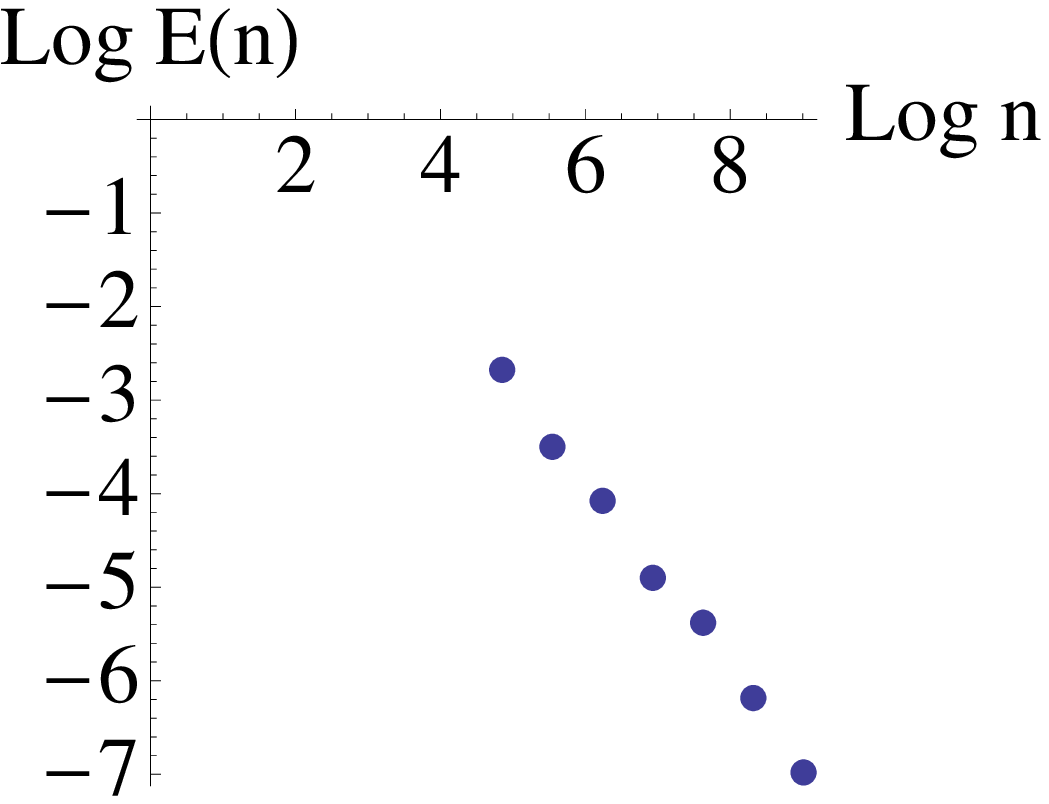}
\caption{Computation of a Dirac impulse (plotted as a blue vertical half line) traversing a smooth ramp $\zeta$, depicted in red.  Left: an initially right moving unit pulse at left of the ramp.  Middle: after the pulse has traversed the ramp, a smooth reflected waveform travels to the left and a mixed waveform, consisting of a Dirac pulse followed by a regular function, travels to the right.  In order to match the scale of the ramp, the amplitude of the waveforms has been scaled up by a factor of 32.  Right: a log-log plot of relative rms error $E(n)$ of the plotted waveform for $n=2^p$, $p=7,8,9,10,11,12,13$.  The slope is very close to $-1$, corresponding to $E(n)\cong 7.4/n$. Computation time is 6 seconds for $p=10$, at which point the plot is indistinguishable from that of the exact result.
}\label{fig-ramp}
}
\end{figure}

Combining Theorems~\ref{thm-two} and \ref{thm-three} one can exploit a separation of scales phenomenon to compute both the regular and singular parts of the solution $u(x,t)$ to (\ref{wave}) in the case of Dirac function initial data.  The solutions to the approximations $u^P$ are Dirac combs whose components have the form $\gamma\delta(x-x_j-t_k)$.  If $(x_j,t_k)$ is a common grid point for a succession of partitions $P$, then the rescaled coefficient $\gamma/\Delta$ diverges if $(x_j,t_k)$ belongs to the singular support of $u$; otherwise $\gamma/\Delta$ converges to the regular value of $u(x_j,t_k)$.\footnote{For technical reasons, one in fact computes the solution at two successive spatial grid points and divides by $2\Delta$.}  On the other hand, the unscaled coefficient $\gamma$ goes to zero at regular points and stabilizes at singular points. See Figure~\ref{fig-ramp}.

\begin{figure}[p]
\hspace*{\fill}
\parbox{6in}{
\includegraphics[width=2in]{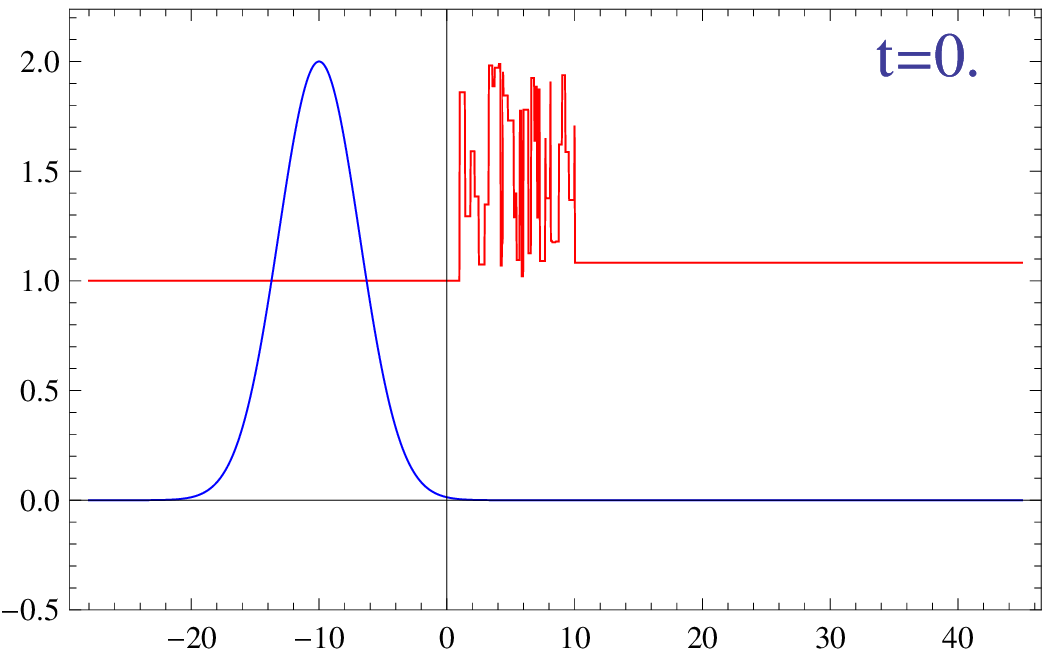}\includegraphics[width=2in]{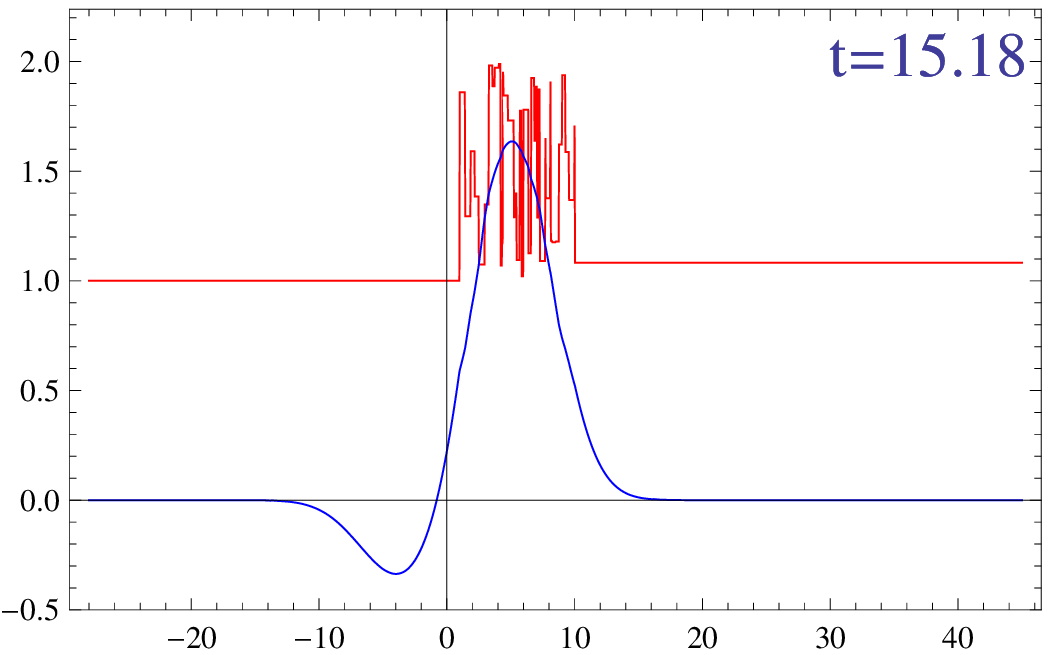}\includegraphics[width=2in]{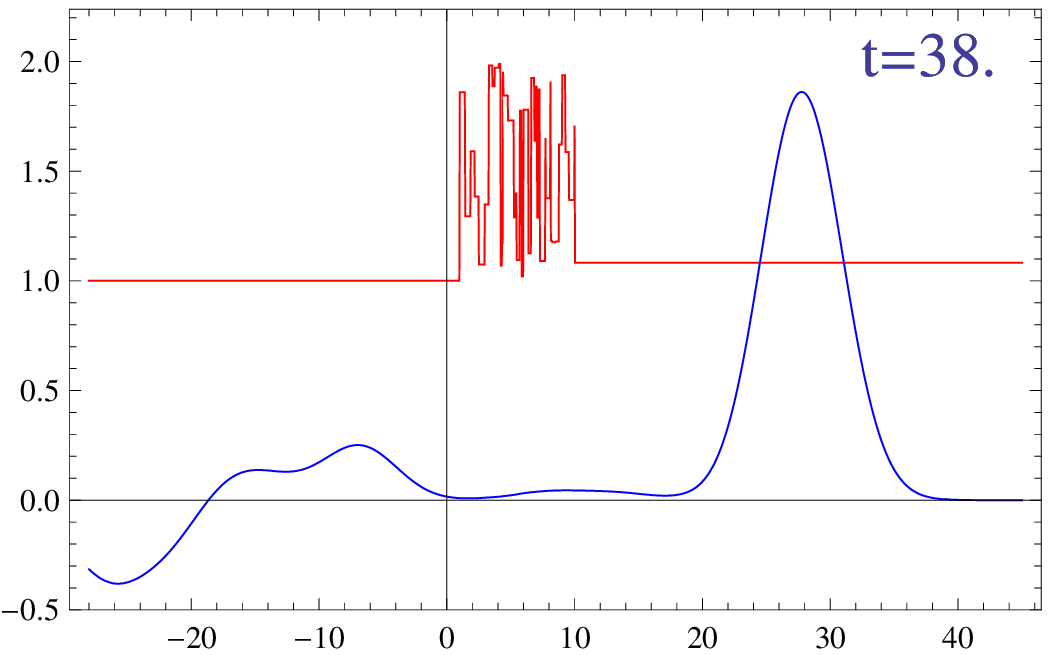}\\
\includegraphics[width=2in]{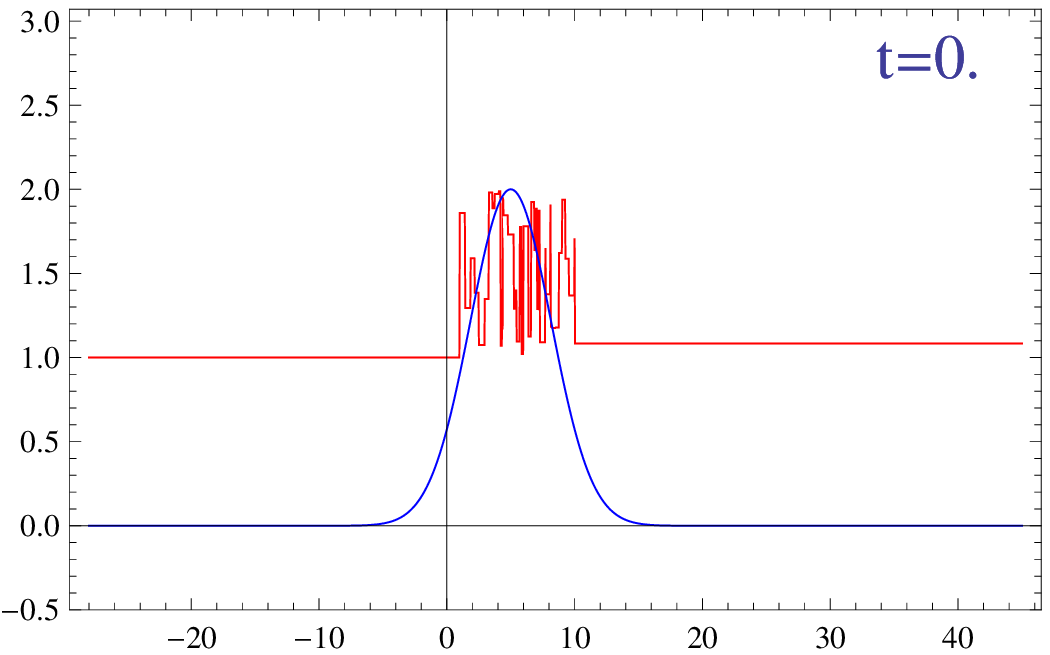}\includegraphics[width=2in]{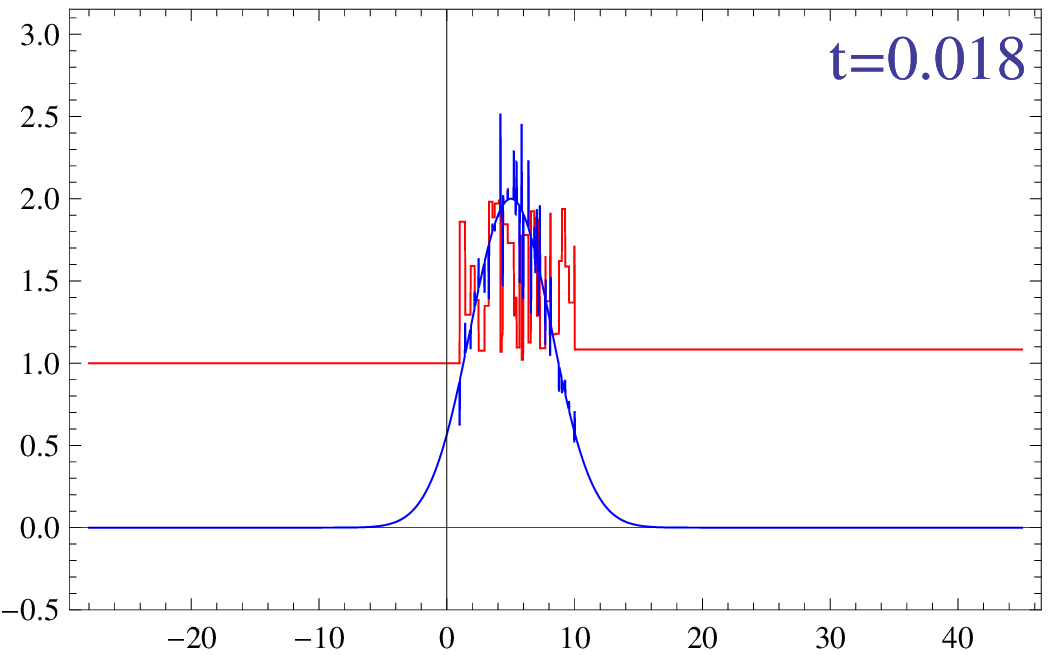}\includegraphics[width=2in]{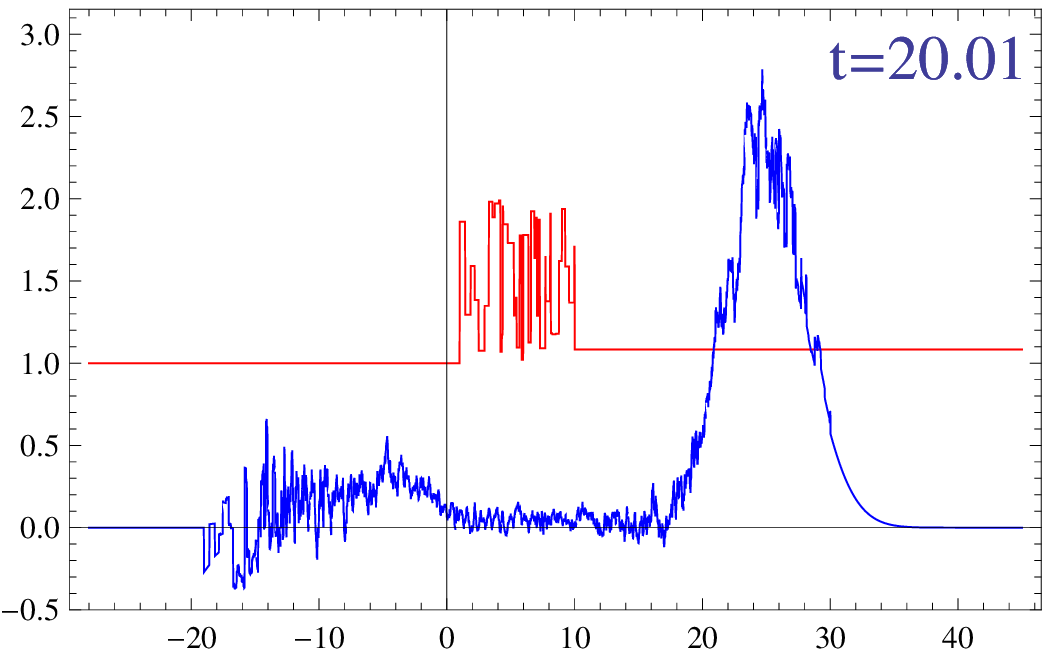}
\caption{The case where $\zeta$, depicted in red, is a highly oscillatory step function.  Top row:  an initially smooth wave form (plotted in blue) traverses the oscillatory zone from left to right, resulting in continuous transmitted and reflected wave forms.  Bottom row: a smooth, right-travelling wave form initially within the oscillatory region spontaneously forms discontinuities, which persist in the resulting right and left-moving wave forms.  This surprising phenomenon may be verified analytically, the simplest case being where $\zeta$ has just a single jump point.}\label{fig-oscillatory}
}
\end{figure}

Algorithm~\ref{alg-main} also handles the contrasting scenario in which initial data is smooth and $\zeta$ is a highly oscillatory step function. If the jump points of $\zeta$ are equally spaced, one can choose $P$ such that $\zeta^P=\zeta$, whereby the algorithm is exact up to roundoff error (by Theorem~\ref{thm-two}).  Comparison with the (much slower) exact methods of \cite{Gi:SIAP2014} shows that a good approximation is also obtained  if the jump points of $\zeta$ are not equally spaced, provided the underlying grid is sufficiently fine.  Indeed, Algorithm~\ref{alg-main} reveals a remarkable qualitative phenomenon whereby smooth initial data can lead to discontinuities, as illustrated in Figure~\ref{fig-oscillatory}. 

The gist of Theorems~\ref{thm-one}, \ref{thm-two} and \ref{thm-three} is that Algorithm~\ref{alg-main} is guaranteed to approximate the solution to the original equation (\ref{wave}) in a wide range of cases, including with distributional initial data. Up to roundoff error, the approximation is exact for piecewise constant $\zeta$, provided its jump points are included in those of $\zeta^P$, irrespective of whether the initial data are smooth or singular. The general requirement that $\zeta$ be a uniform limit of step functions (i.e., a regulated function) is very weak, allowing $\zeta$ to have essentially arbitrary structure on $[x_-,x_+]$ (see \cite[Ch.~VII]{Di:1960}).  Thus Algorithm~\ref{alg-main} differs essentially from existing methods in its interpretation---initial data may be non-smooth or purely distributional, and the coefficient $\zeta$ may have discontinuities.  

It does not appear that Algorithm~\ref{alg-main} can be generalized to higher dimensions in a way that retains the exactness results of Theorem~\ref{thm-two}---in this sense the method is purely one dimensional.  Within the context of one-dimensional evolutionary systems governed by (\ref{wave}) however, the algorithm marks a qualitative improvement over existing methods in terms of both speed and generality.  
In summary, Algorithm~\ref{alg-main} comprises an easily-implemented, versatile and accurate computational tool applicable to various imaging modalities and to modelling of various biological or built structures.  


\appendix

\section{Computational data}

The present section gives further details on Figures~\ref{fig-ramp} and \ref{fig-oscillatory} to facilitate replication of the results.  The smooth ramp in Figure~\ref{fig-ramp} is given by the formula
\[
\zeta(x)=\left\{\begin{array}{cc}
1&\mbox{ if }x\leq 1\\
\rule[-10pt]{0pt}{28pt}2+\tanh\left(\frac{8(x-5)}{16-(x-5)^2}\right)&\mbox{ if }1<x<9\\
3&\mbox{ if }9\leq x
\end{array}\right..
\]
The source is a right moving unit Dirac function initially centered at $x=-5$.  

In the case of Figure~\ref{fig-oscillatory}, $\zeta$ is a step function with 40 randomly generated jump points lying  in the interval $[1,10]$. $\zeta(x)=1$ for $x<1$ and $\zeta(x)=2/3$ for $x>10$.  The source wave is a right-travelling gaussian initially of the form 
\[
\alpha(x)=2e^{-.05(x+10)^2}
\]
in the first example; the source is shifted 15 units to the right in the second example.  
The qualitative nature of the propagating wavefield does not depend on the detailed structure of the step function.

\end{document}